\documentclass[12pt]{article}
\usepackage{amssymb,amsmath}

\author{Abdallah Assi\thanks{Universit\'e d'Angers, Math\'ematiques,
49045 Angers ceded 01, France, e-mail:assi@univ-angers.fr \break Visiting address: American
University of Beirut, Department of Mathematics, Beirut 1107 2020,
Lebanon}}
\title{Rational curves with one place at infinity
}
\date{\mbox{}}

\newtheorem{teorema}{Theorem}[section]

\newtheorem{proposicion}[teorema]{Proposition}
\newtheorem{lema}[teorema]{Lemma}

\newtheorem{nota}[teorema]{Remark}
\newenvironment{demostracion}[1]{\paragraph{\sl Proof#1}}{}
\newenvironment{demostracione}[1]{\paragraph{\sl Proof of Theorem 3.1.#1}}{}

\newcommand{\qed}{\mbox{$\Box$}}

\begin{document}
\maketitle

\noindent{\bf Abstract:} Let ${\mathbb K}$ be an algebraically closed field of characteristic zero. Given a polynomial $f(x,y)\in{\mathbb K}[x,y]$  with one place at infinity, we prove that either $f$ is equivalent to a coordinate, or the family $(f_{\lambda})_{\lambda\in\mathbb{K}}$ has at most two rational elements. When $(f_{\lambda})_{\lambda\in\mathbb{K}}$ has two rational elements, we give a description of the singularities of these two elements.

\medskip

\section{Introduction and notations }
 
\medskip

\noindent Let ${\mathbb K}$ be an algebraically closed field of characteristic zero, and let $f=y^n+
a_1(x)y^{n-1}+\ldots+a_n(x)$ be a monic reduced polynomial of ${\mathbb
  K}[x][y]$. For all $\lambda\in {\mathbb K}$, we set $f_{\lambda}=f-\lambda$. Hence we get a family of polynomials $(f_{\lambda})_{\lambda\in {\mathbb K}}$. We shall suppose that $f_{\lambda}$ is a reduced polynomial for all $\lambda\in {\mathbb K}$.  Let $g$ be a nonzero polynomial of ${\mathbb K}[x][y]$.  We define the intersection
multiplicity of $f$ with $g$, denoted int$(f,g)$, to be the rank of the ${\mathbb K}$-vector space $\displaystyle{{{{{\mathbb
    K}[x][y]}}\over {(f,g)}}}$. Note that int$(f,g)$ is also the
$x$-degree of the $y$-resultant of $f$ and $g$.  Let $p=(a,b)\in
V(f)\cap V(g)$, where $V$ denotes the set of zeros in ${\mathbb K}^2$. By setting $\bar{x}=x-a,\bar{y}=y-b$, we may   assume  that $p=(0,0)$. We define the intersection multiplicity of $f$ with $g$ at $p$, denoted int$_p(f,g)$,  to be the rank of the ${\mathbb K}$-vector space $\displaystyle{{{{{\mathbb
    K}[[x,y]]}}\over {(f,g)}}}$.  Note that int$(f,g)=\sum_{p\in  V(f)\cap V(g)}{\rm int}_p(f,g)$. We define the local Milnor number of $f$ at $p$, denoted $\mu_p(f)$,  to be the intersection multiplicity int$_p(f_x,f_y)$, where $f_x$ (resp. $f_y$) denotes the $x$-derivative (resp. the $y$-derivative) of $f$. We set $\mu(f)=\sum_{p\in V(f)}\mu_p(f)$ and $\mu={\rm int}(f_x,f_y)$ and we recall that $\mu=\sum_{\lambda\in{\mathbb K}}\mu(f_{\lambda})=\sum_{\lambda\in{\mathbb K}}\sum_{p\in V(f_{\lambda})}\mu_p(f_{\lambda})$. Let $q$ be a point in $V(f)$ and assume, after possibly a change of
variables that $q=(0,0)$. The number of places of $f$ at $q$, denoted
$r_{q}$, is defined to be the number of irreducible components of
$f$ in ${\mathbb K}[[x,y]]$.

\noindent Assume, after possibly a change of variables, that
deg$_xa_i(x)<n-i$ for all $i=1,\ldots,n$ (where deg$_x$ denotes the
$x$-degree). In particular $f$ has one point at infinity defined by
$y=0$. Let $h_f(x,y,u)=u^n\displaystyle{f({x\over u},{y\over
    u})}$. The local equation of $f$ at infinity  is nothing but
$F(y,u)=h_f(1,y,u)\in{\mathbb K}[[u]][y]$. We define the Milnor number
of $f$ at infinity, denoted $\mu_{\infty}$, to be the rank of the
${\mathbb K}$-vector space $\displaystyle{{{\mathbb K}[[u]][y]}\over
  {(F_u,F_y)}}$. We define the number of places at infinity of $f$, denoted
$r_{\infty}$, to be the number of irreducible components of
$F(y,u)$ in ${\mathbb K}[[u]][y]$.

\section {Curves with one place at infinity}

\medskip

\noindent Let the notations be as in Section 1., in particular $f=y^n+a_1(x)y^{n-1}+\ldots+a_n(x)$ is a monic reduced 
polynomial of ${\mathbb K}[x,y]$. Let
$R(x,\lambda)=P_0(\lambda)x^i+\ldots+P_i(\lambda)$ be the $y$-resultant
of $f_\lambda,f_y$. We say that $(f_{\lambda})_{\lambda\in{\mathbb K}}$ is $d$-regular (discriminant-regular)
if $P_0(\lambda)\in\mathbb{K}^*$. Note that $(f_{\lambda})_{\lambda\in{\mathbb K}}$
is $d$-regular if and only if int$(f_{\lambda},f_y)=i$ for all
$\lambda\in\mathbb{K}$. Suppose that $(f_{\lambda})_{\lambda\in{\mathbb K}}$ is not
$d$-regular, and let $\lambda_1,\ldots,\lambda_s$ be the set of roots of
$P_0(\lambda)$. We set $I(f)=\lbrace
\lambda_1,\ldots,\lambda_s\rbrace$, and we call $I(f)$ the set of
$d$-irregular values of $(f_{\lambda})_{\lambda\in{\mathbb K}}$. Let 
$A_f=\sum_{k=1}^s(i-{\rm int}(f-\lambda_k,f_y))$. For all $\lambda\in \mathbb{K}-I(f)$, we have int$(f_{\lambda},f_y)=\mu+n-1+A_f$,  where $\mu={\rm int}(f_x,f_y)$ (see [5]).

\medskip

\noindent  Note that $A_f=\sum_{\lambda\in\mathbb{K}}(i-{\rm int}(f_{\lambda},f_y))$, in particular $(f_{\lambda})_{\lambda\in{\mathbb K}}$ is $d$-regular if and only if $A_f=0$. On the other hand, given $a\in\mathbb{K}$, if int$(f_{a},f_y)=\mu+n-1$, then either $(f_{\lambda})_{\lambda\in{\mathbb K}}$ is $d$-regular or $I(f)=\lbrace a\rbrace$.

\medskip
\noindent Assume that deg$_xa_k(x) < k$ for all $k=1,\ldots,n$, in such a way that $y=0$ is the only point at infinity of $f$.

\begin{proposicion}{\rm (see [2] and [3]) Let the notations be as above and assume that  $f$ has one place at infinity,
i.e. the projective curve defined by the homogeneous equation
$h_f(x,y,u)=\displaystyle{f({x\over u},{y\over u})u^n}$
is analytically irreducible at the point at infinity $(1:0:0)$.  We have the following

\begin{itemize}

\item For all $\lambda\in{\mathbb K}, f-\lambda$ has one place at infinity.

\item The family $(f_{\lambda})_{\lambda\in{\mathbb K}}$ is $d$-regular. In particular, int$(f_{\lambda},f_y)=\mu+n-1$ for all $\lambda\in{\mathbb K}$.

\item If $\mu=0$, then deg$_xa_n(x)$ divides $n$ and there exists an automorphism $\sigma$ of ${\mathbb K}^2$ such that $\sigma(f)$ is
a coordinate of ${\mathbb K}^2$.
\end{itemize}
}
\end{proposicion}

\noindent Let the notations be as above. If $\delta_p$ (resp. $\delta_{\infty}$) denotes the order of the
conductor of $f$ at $p\in V(f)$ (resp. at the point at infinity), then
$2\delta_p=\mu_p+r_p-1$
(resp. $2\delta_{\infty}=\mu_{\infty}+r_{\infty}-1$) (see [7]).  Assume that $f$ is an irreducible polynomial, and let $g(f)$
be the genus of the normalized curve of $V(f)$. By the genus formula
we have: 

$$
2g(f)+(\sum_{p\in
  V(f)}2\delta_p)+2\delta_{\infty}=(n-1)(n-2).
$$

\noindent Now  int$(f,f_y)=\mu+n-1+A(f)$, where
$A(f)$ is a nonnegative integer and $A(f)=0$ if and only if 
$(f_{\lambda})_{\lambda\in {\mathbb K}}$ has at most one $d$-irregular value at
infinity. On the other hand, the local intersection multiplicity of $f$ with $f_y$
 at the point at infinity is $\mu_{\infty}+n-1$. In particular
$\mu+\mu_{\infty}=(n-1)(n-2)$,
consequently, if  $\mu(f)=\sum_{p\in V(f)}\mu_p$, and
$\overline{\mu}(f)=\mu-\mu(f)$, then 

$$
2g(f)+(\sum_{p\in
  V(f)}2\delta_p)+2\delta_{\infty}=\mu(f)+\overline{\mu}(f)+\mu_{\infty}+A(f).
$$

\noindent We finally get: 

$$
(**)\quad 2g(f)+\sum_{p\in
  V(f)}(r_p-1)+r_{\infty}-1=\overline{\mu}(f)+A(f)
$$

\noindent in particular $g(f)=\sum_{p\in
  V(f)}(r_p-1)+r_{\infty}-1=0$ if and only if
$A(f)=\overline{\mu}(f)=0$. Roughly speaking, $f$ is a rational
unibranch curve (at infinity as well as at finite distance) if and
only if the pencil $(f_{\lambda})_{\lambda\in {\mathbb K}}$ has at most one $d$-irregular value at infinity and for all $\lambda\not=0, f_{\lambda}$ is a smooth curve. Under these hypotheses, Lin-Zaidenberg Theorem implies  that $f$ is equivalent to a quasihomogeneous curve $Y^a-X^b$ with gcd$(a,b)=1$ (see [6]). Note that these hypotheses are satisfied when $r_{\infty}-1=0=\mu$. Hence we get the third assertion of Proposition 2.1. since in this case, min$(a,b)=1$ and $f$ is equivalent to a coordinate

\section{Rational one place curves}

\noindent  Let $f=y^n+a_1(x)y^{n-1}+\ldots+a_n(x)$ be a polynomial of ${\mathbb K}[x,y]$ and let the notations be as in Sections 1 and 2. Assume that $f$ has one place at infinity, i.e. $r_{\infty}=1$. If $f$ is rational, then it follows from the equality (**) of Section 2 that $\sum_{p\in
  V(f)}(r_p-1)=\overline{\mu}(f)$. We shall prove the following:


\begin{teorema}{\rm  Assume that $f$ has one place at infinity and let 
    $(f_{\lambda})_{\lambda\in{\mathbb K}}$ be the pencil of curves defined by
    $f$ . If $f$ is rational, then exactly one of the following holds:

i) For all $\lambda\in\mathbb{K}$, $f_{\lambda}$ is rational, and $\sigma(f)$  is a coordinate of $\mathbb{K}^2$ for some automorphism $\sigma$ of $\mathbb{K}^2$.

ii) The polynomial $f-\lambda$ is rational for at most one
$\lambda_1\not=0$, i.e. the pencil  $(f_{\lambda})_{\lambda\in{\mathbb K}}$ has at
most two rational elements.

}
\end{teorema}

\noindent We shall  prove first the following Lemma:

\begin{lema}{\rm Let $H=y^N+a_1(x)y^{N-1}+\ldots+a_N(x)$ be a non zero reduced polynomial of ${\mathbb
      K}[[x]][y]$, and let $H=H_1\ldots H_r$ be the decomposition of
    $H$ into irreducible components of  ${\mathbb K}[[x]][y]$. Let
    $\mu_{(0,0)}$ denotes the Milnor number of $H$ at $(0,0)$
    (i.e. $\mu_{(0,0)}$ is the rank of the ${\mathbb K}$-vector space
    $\displaystyle{{{\mathbb K}[[x]][y]}\over {(H_x,H_y)}})$. We have
    the following:

i) $\mu_{(0,0)}\geq r-1$. 

ii) If $r\geq 3$, then $\mu_{(0,0)}>r-1$.

iii) If $r=2$ and $\mu_{(0,0)}= r-1=1$, then $(H_1,H_2)$ is a local system of coordinates at $(0,0)$.
}
\end{lema}

\begin{demostracion}{.} We have int$_{(0,0)}(H,H_y)=\mu_{(0,0)}+N-1$,
  but
$$
 {\rm int}_{(0,0)}(H,H_y)=\sum_{i=1}^r{\rm
   int}(H_{i},H_{i_y})+2\sum_{i\not= j}{\rm int}_{(0,0)}(H_i,H_j)
$$

$$
=\sum_{i=1}^r{\rm
   int}[(H_{i_x},H_{i_y})+{\rm deg}_yH_i-1]+2\sum_{i\not= j}{\rm
   int}_{(0,0)}(H_i,H_j)
$$

\noindent hence

$$
\mu_{(0,0)}+N-1=(\sum_{i=1}^r{\rm
   int}(H_{i_x},H_{i_y}))+N-r+2\sum_{i\not= j}{\rm
   int}_{(0,0)}(H_i,H_j).
$$

\noindent Finally we have $\mu_{(0,0)}=(\sum_{i=1}^r{\rm
   int}(H_{i_x},H_{i_y}))-r+1+2\sum_{i\not= j}{\rm
   int}_{(0,0)}(H_i,H_j)$ . Now  for all $1\leq i\leq r$,  ${\rm
   int}_{(0,0)}(H_{i_x},H_{i_y})\geq 0$ and $\sum_{i\not= j}{\rm
   int}_{(0,0)}(H_i,H_j)\geq C_2^r=\displaystyle{{r(r-1)}\over 2}$,
 hence  $\mu_{(0,0)}\geq r(r-1)-(r-1)=(r-1)^2$ and i), ii) follow
 immediately. Assume that $r=2$. If $\mu_{(0,0)}=r-1$, then  ${\rm
   int}_{(0,0)}(H_{1_x},H_{1_y})= {\rm
   int}_{(0,0)}(H_{2_x},H_{2_y})=0$ and int$_{(0,0)}(H_1,H_2)=1$. This implies iii)$\blacksquare$
\end{demostracion}

\begin{demostracione}{.} {\rm  If $\mu(f)=0$, then $\mu=0$ and by
  Proposition 2.1.,  $\sigma(f)$ is a coordinate of $\mathbb{K}^2$ for some
  automorphism $\sigma$ of $\mathbb{K}^2$. Assume that $\mu(f)>0$ and
  let $p_1,\ldots,p_s$ be the set of singular points of $V(f)$. Let
  $r_i$ denotes the number of places of $f$ at $p_i$ for all $1\leq
  i\leq s$. By
  Lemma 3.2., for
  all $1\leq i\leq s$, $\mu_{p_i}\geq r_i-1$, on the other hand,
  equality (**) of Section 2 implies that
  $\sum_{i=1}^s(\mu_{p_i}+r_i-1)=\mu$, in particular $\mu \leq
  \sum_{i=1}^s2\mu_{p_i}=2\mu(f)$, hence $\mu(f)\geq
  \displaystyle{\mu\over 2}$. If $f_{\lambda_1}$ is rational for some
  $\lambda_1\not=0$, then the same argument as above implies that  $\mu(f_{\lambda_1})\geq
  \displaystyle{\mu\over 2}$. This is possible only for at most one
  $\lambda_1\not= 0$, hence  ii) follows immediately.
$\blacksquare$
}
\end{demostracione}

\noindent The following proposition characterizes the case where the pencil $(f_{\lambda})_{\lambda\in{\mathbb K}}$ 
has exactly two rational elements.

\begin{proposicion}{\rm Let the notations be as in Theorem 3.1.   and assume that the pencil  $(f_{\lambda})_{\lambda\in{\mathbb K}}$ has exactly two
rational elements $f$ and $f_{\lambda_1}$. We have $\mu(f)=\mu(f_{\lambda_1})=
  \displaystyle{\mu\over 2}$, furthermore, given a singular point $p$
  of $V(f)$ (resp. $V(f_{\lambda_1})$), $f$ (resp. $f_{\lambda_1}$)
  has two places at $p$ and $\mu_p(f)=1$ (resp. $\mu_p(f_{\lambda_1})=1$). In particular, 
$f$ (resp. $f_{\lambda_1}$) has exactly $\displaystyle{\mu\over 2}$ singular points.
}
\end{proposicion}

\begin{demostracion}{.}  It follows from the proof of Theorem 3.1. that $\mu(f)\geq \displaystyle{\mu\over 2}$ and that 
$\mu(f_{\lambda_1})\geq \displaystyle{\mu\over 2}$. Clearly this holds only if $\mu(f)=\mu(f_{\lambda_1})=
  \displaystyle{\mu\over 2}$. Let $p$ be a singular point of $V(f)$. We have $\mu_p=r_p-1$, hence, by Lemma 3.2. ii), 
$r_p\leq 2$. But $\mu_p>0$, hence $r_p=2$ and $\mu_p=1$. This implies that $f$ has $\displaystyle{\mu\over 2}$ 
singular points. Clearly the same holds for $f_{\lambda_1}$.$\blacksquare$
\end{demostracion}

\noindent The results above imply  the following:

\begin{proposicion}{\rm Assume that $f$ has one place at infinity and  let $(f_{\lambda})_{\lambda\in{\mathbb K}}$ be the pencil of polynomials defined by
    $f$. Assume that $f$ is a rational polynomial and that $\mu(f) >0$.  Let $p_1,\ldots,p_s$ be the set of singular points of $f$. We have the following

i) If $r_{p_i}=1$ (resp. $r_{p_i}\geq 3$)  for some $1\leq i\leq s$, then $f$ is the only rational point of the pencil $(f_{\lambda})_{\lambda}$ .

ii) If $r_{p_i}=2$ for all $1\leq i\leq s$ but $s\not= \displaystyle{\mu\over 2}$,  then $f$ is the only rational element of the pencil $(f_{\lambda})_{\lambda}$ .

}
\end{proposicion}

\begin{demostracion}{.} This is an immediate application of Theorem 3.1. and Proposition 3.3.$\blacksquare$
\end{demostracion}

\begin{proposicion}{\rm Let $f\not=g$ be two monic polynomials of ${\mathbb K}[x][y]$ and assume that $f,g$ are parametrized by polynomials of ${\mathbb K}[t]$. Under these hypotheses, exactly one of the
    following conditions holds:

i) $f=g+\lambda_1$ for some $\lambda_1\in\mathbb{K}^*$, and $f$ is
equivalent to a coordinate, i.e. $\sigma(f)$ is a coordinate of
$\mathbb{K}^2$ for some automorphism $\sigma$ of $\mathbb{K}^2$.

ii) $f=g+\lambda_1$ for some $\lambda_1\in\mathbb{K}^*$, $\mu(f)=\mu(g)=
  \displaystyle{{{\rm int}(f_x,f_y)}\over 2} >0$, and $f$ (resp. $g$) has $\displaystyle{{{\rm int}(f_x,f_y)}\over 2} $
singular points with two places at each of them.

iii)  int$(f,g) > 0$, i.e. $f,g$ meet in a least one point of
$\mathbb{K}^2$.
}
\end{proposicion}

\begin{demostracion}{.}  The polynomial $f$ (resp. $g$) has one place at infinity.  If int$(f,g)=0$, then 
$f=ag+\lambda_1, a,\lambda_1\in{\mathbb K}^*$. Since $f$ and $g$ are monic, then $a=1$. Hence 
$g$ and $g+\lambda_1$ are two rational elements of the pencil $(f_{\lambda})_{\lambda\in{\mathbb K}}$. Now apply Theorem 3.1.
and Proposition 3.3.$\blacksquare$
\end{demostracion} 

\begin{nota}{\rm Let $(x(t),y(t))=(t^3-3t,t^2-2)$ and
    $(X(s),Y(s))=(s^3+3s,s^2+2)$, and let $f(x,y)={\rm
      res}_t(x-x(t),y-y(t))$ (resp. $g(x,y)={\rm
      res}_s(x-X(s),y-Y(s))$). We have $(x(t)-X(s),y(t)-Y(s))={\mathbb
      K}[t,s]$, hence int$(f,g)=0$. In fact,
    $$
f(x,y)=y^3-x^2-3y+2=-x^2+(y+2)(y-1)^2
$$ 

\noindent and
    
$$
g(x,y)=y^3-x^2-3y-2=
-x^2+(y-2)(y+1)^2,
$$

\noindent hence $f=g+4$. The genus of a generic element of the family
$(f_{\lambda})_{\lambda}$ is $1$, and $f,f-4$ are the two rational
elements of this family. Note that $\mu=2$ and $\mu(f)=\mu(f-4)=1$.
This example shows that the bound of Theorem 3.1. is sharp. 
}
\end{nota}

\begin{nota}{ \rm Let $(f_{\lambda})_{\lambda\in{\mathbb K}}$ be a pencil of
    polynomials of ${\mathbb K}[x,y]$ and   assume that $f-\lambda$ is
    irreducible for all $\lambda\in {\mathbb K}$. If the generic
    element of the pencil is rational, then for all
    $\lambda\in{\mathbb K}$, $f-\lambda$ is rational and
    irreducible. In this case, by  [8], $f$ has one place at infinity
    and $\sigma(f)$ is a coordinate of ${\mathbb K}^2$ for some
    automorphism $\sigma$ of ${\mathbb K}$. Assume that the genus of the generic element of the pencil $(f_{\lambda})_{\lambda\in{\mathbb K}}$ is greater than or equal to one. Similarly to the case of curves with one place at infinity, it is natural to address the following question:







\noindent {\bf Question}: Is there an integer $c\in{\mathbb N}$ such
that, given a pencil of irreducible polynomials
$(f_{\lambda})_{\lambda\in{\mathbb K}}$, if  $\mu+A_f >0$, then the number of
rational elements in the pencil is bounded by $c$?
}


\end{nota}
\medskip

\noindent Acknowledgments: The author would like to thank the referee for the valuable and helpful comments.

\noindent [1] S.S. Abhyankar.- On the semigroup of a meromorphic
curve, Part 1, in Proceedings, International Symposium on Algebraic
Geometry, Kyoto (1977), 240-414.

\noindent [2] S.S. Abhyankar and T.T. Moh.- Newton-Puiseux expansion
and generalized Tschirnhausen transformation, Crelle Journal, 260 (1973)
47-83.

\noindent [3] S.S. Abhyankar and T.T. Moh.- Newton-Puiseux expansion, 
and generalized Tschirnhausen transformation II, Crelle Journal,
261 (1973), 29-54.

\noindent [4] A. Assi.- Sur l'intersection des courbes m\'eromorphes,
C. R. Acad. Sci. Paris S\'er. I Math. 329 (1999), n$^0$7, 625-628.

\noindent [5] A. Assi.- Meromorphic plane curves,
Math. Z. 230 (1999), n$^0$1, 16-183.


\noindent [6] V. Lin and M. Zaidenberg.- An irreducible simply connected algebraic curve in ${\mathbb C}^2$ is
equivalent to a quasihomogeneous curve, Dokl. Akad. Nauk SSSR, 271 (1983), n$^0$5, 1048-1052.

\noindent [7] J. Milnor.- Singular points of complex hypersurfaces,
Ann. of Math. Studies, 61, Princeton, Univ. Press., Princeton, NJ,
1968.

\noindent [8] W. Neumann and P. Norbury.- Nontrivial rational polynomials in two variables have reducible fibres,
Bull. Austral. Math. Soc., Vol 58 (1998), 501-503.

\noindent [9] O. Zariski.- Le probleme des modules pour les
branches planes, Lectures at Centre de Math\'ematiques, Ecole
Polytechnique, Notes by F. Kmety and M. Merle, 1973.

\end{document}